\documentclass[12pt]{article}

\usepackage{amsmath}
\usepackage{amssymb}
\usepackage{latexsym}
\usepackage{amsfonts}
\usepackage{fullpage}
\usepackage{verbatim}
\usepackage[utf8]{inputenc}

\newtheorem{thm}{Theorem}

\begin{document}

\centerline {\hfil Note on Artin's Conjecture on Primitive Roots \hfil }

\bigskip

\centerline { \hfil Sankar Sitaraman \hfil }

\centerline { \hfil Howard University, Washington, DC \hfil }

\bigskip
\centerline { \hfil \it To  Gerald, Ralph, Sundar and Susan \hfil } 
\bigskip

\centerline {\bf Abstract}

\smallskip

 E. Artin conjectured that any integer
$a > 1$ which is not a perfect square is a primitive root modulo $p$ for infinitely many primes $ p.$ Let $f_a(p)$ be the  multiplicative order of the non-square integer $a$ modulo the prime $p.$ M. R. Murty and S. Srinivasan \cite{Murty-Srinivasan} showed that if $\displaystyle \sum_{p < x} \frac 1 {f_a(p)} = O(x^{1/4})$ then  Artin's conjecture is true for $a.$ We relate the Murty-Srinivasan condition to sums involving the cyclotomic periods from the subfields of $\mathbb Q(e^{2\pi i /p})$ corresponding to the subgroups $<a> \subseteq \mathbb F_p^*.$ 

\bigskip

\bigskip

\centerline {\bf 1. Introduction}

\medskip

\bigskip

Let $a>1$ be an integer
 which is not a perfect square and let $f_a(p)$ be the  multiplicative order of  $a$ modulo a given prime $p.$  For such $a$  the   Murty-Srinivasan \cite{Murty-Srinivasan} condition $\displaystyle \sum_{p < x} \frac 1 {f_a(p)} = O(x^{1/4})$  applies, and one could try to evaluate this sum by looking at the average value of $ \omega(\Phi_j(a))$ where $\Phi_j(x)$ is the $j-$th cyclotomic polynomial and $\omega(n)$ is the number of distinct primes dividing $n.$ We were led to this idea by our thinking of $f_a(p)$ as the lengths of the orbits of the points on the unit circle $S^1$ under the ``Frobenius" map $\sigma_a: x \to x^a.$ This meant looking at the torsion of the group of roots of unity. Note that for a given prime $p,$  $f_a(p)$ is the order of this map acting on the multiplicative subgroup generated by $\zeta_p = e^{\frac {2\pi i }p}.$ It is well known that if $f_a(p) = j$ then $p \ | \ \Phi_j(a)$ where $\Phi_j(x)$ is the $j-$th cyclotomic polynomial. See, for instance, \cite{Washington}, Chapter 2. 

 Since the $\Phi_j(a)$ divide $a^j-1$ and we are only concerned with $p < x$ we actually use  $\sum \omega_x(a^j-1)$ where $\omega_x(m)$ for integral $m$ is the number of primes upto $x$ that divide it.  The argument to approach the Murty-Srinivasan condition using the average value of $ \omega(\Phi_j(a))$  is also  in Murty-Wong \cite{Murty-Wong}, Section 5, albeit implicitly. There an inequality for a sum that contains $S(x)= \displaystyle \sum_{p < x} \frac 1 {f_a(p)} $ is used to give heuristics for the largest prime factor of $a^j-1, j \le x.$ We show how the sum $\displaystyle \sum_{j \le x} {\omega_x(a^j-1)} $ can be related to cyclotomic periods after expressing it in terms of  finite Fourier  series (exponential sums) and give some bounds for these sums, from which one could easily get bounds for  $\displaystyle \sum_{p < x} \frac 1 {f_a(p)}.$

 \vfill\eject

\noindent {\bf Statement of Main Result: Theorem 1}

\medskip

In what follows $f_a(p)$ will be sometimes written simply as $f_p$ where the context is clear, and 
$\displaystyle e_p $ will denote $ \displaystyle \frac {p-1}{f_p}.$ 

$g$ will denote a primitive generator of $\mathbb F_p^*$ and so $a \equiv  g^{e_p} \pmod{ p}.$
$g_i $ will denote the least element of $\{1,2,...,p-1\}$ such that $g_i \equiv g^i \pmod{ p}.$

The cyclotomic periods for $i = 1,2,..., e_p,$ denoted $\eta_{i,p},$ are defined by
$$ \eta_{i,p} = \sum_{j=1}^{f_p} \zeta_p^{g^i g^{je_p}} = \sum_{j=1}^{f_p} \zeta_p^{g^i a^j} = \sum_{j=1}^{f_p} \zeta_p^{g_i a^j}.$$  

They are contained in the subfield of $\mathbb Q(\zeta_p)$ corresponding to the subgroup $<a>$ generated by $a$ in $\mathbb F_p^*.$
$H_p$ will denote this subgroup. 

  \bigskip

 \begin{thm} For $a,$  $f_a(p)$ and $\eta_{i,p}$ corresponding to $H_p$ defined as above,  we get the following bounds for any $u \le x$ where $x$ is a fixed positive integer:
 
  $$A: \quad \sum_{j \le u} {\omega_x(a^j-1)}  \le    O(u\log \log x)+ u \sum_{p \le x} \frac {\displaystyle \sum_{ i = 1} ^{\frac{p-1}{f_a(p)}} |\eta_{i,p}|}p.$$

 $$ B: \quad \quad \sum_{j \le u} {\omega_x(a^j-1)}  \le O(u \log \log x) + 2\pi u \sum_{p \le x} \frac { \left(\sum_{k=1}^{p-2} \left | \sum_{m=1}^k \eta_{h_m,p} \right | \right) }  {p^2} $$
 where $\displaystyle \eta_{h_m,p}= \sum_{j=1}^{f_a(p)} \zeta_p^{h_ma^j}$ with $h_m \equiv mf_a(p) \pmod{p}, m \in [1,p-1].$
 
 \end{thm}

Much progress has been made on Artin's conjecture by Hooley \cite{Hooley}, Gupta-Murty \cite{Gupta-Murty}, Heath-Brown \cite{Heath-Brown} and others. Heath-Brown proved that Artin's conjecture is true for one of 2, 3 and 5. 

Our initial motivation was to come up with a Dirichlet series that would help in obtaining information about the asymptotic behavior of invariants of cyclotomic fields such as class numbers. We hope to provide such a framework that would enable us to calculate the asymptotics of $f_a(p)$ for a given  $a$ more precisely, as indicated in Sitaraman \cite{Sitaraman}. For instance,  Artin conjectured that $ f_a(p) = p-1$ with a certain non-zero density $C$ (equal to 0.374 approximately) known now as the Artin constant.

 The author is grateful to Niranjan Ramachandran and Larry Washington for their hospitality and conversations during his visit to the University of Maryland, College Park, in 2014-15. He would also like to  thank  Francois Ramaroson and Dani Szpruch for helpful conversations, and the anonymous referee for valuable suggestions to improve the paper.
 
\vfill\eject

\centerline {\bf 2. From the Murty-Srinivasan criterion to Cyclotomic Periods. }

\bigskip

\noindent {\bf 2a. Reduction to estimate of $\sum \omega_x(a^j-1)$}

\medskip

Let $\displaystyle S(x) = \sum_{p < x} \frac 1 {f_a(p)}. $ Let $N(j)$ to be the number of
primes $p < x $ with $a$ of order $ j$ modulo $p.$ (i.e, $f_a(p) = j$). Then we have 
$ \displaystyle S(x) \le  \sum_{j < x} \frac {N(j)} j.$ Since as mentioned before $f_a(p) = j \implies p | \Phi_j(a)$  we get
$$ \displaystyle S(x) = \sum_{j < x} \frac {N(j)} j \le \sum_{j < x} \frac {\omega_x(\Phi_j(a))} j.$$ 

Since $\Phi_j(a)$ is a factor of $a^j-1$ we have $ \displaystyle  \sum_{j < x} \frac {\omega_x(\Phi_j(a))} j \le  \sum_{j < x}\frac {\omega_x(a^j-1)} j .$ Since it is easier to deal with $a^j-1$ we consider the sum $ \displaystyle  \sum_{j < x} \frac {\omega_x(a^j-1)} j .$ We will denote this sum by  $W(x).$
\medskip

Initially we tried to use sieve theoretic techniques to deal with $W(x).$ H. Halberstam \cite{Halberstam} obtained an Erd\"os-Kac type result for $ \displaystyle  \omega(f(m))$ where $f(x) \in \mathbb Z[x]$ is an irreducible polynomial. Namely he showed that 
$ \displaystyle  \sum_{m \le n} \frac {\omega(f(m)) - A(n) }{n \sqrt{A(n)} } $ is asymptotically normally distributed, where $ \displaystyle A(n) =  \sum_{p \le n} \frac {r(p)}p$ with $r(p)$ being the number of solutions to  $f(k) \equiv 0 \ \text {mod}\  p$ with $0 \le k < p.$ Halberstam needs $r(p)$ to be uniformly bounded, and that is possible because the polynomial is fixed. Here we need to average over several polynomials of varying degree.  Even using a later more general version of Halberstam's result proved by Granville-Soundararajan \cite{Granville-Sound} we found it difficult to get a decent estimate for $W(x).$

On the other hand A. T.  Felix \cite{Felix} has shown that on average $\sum_{p<x}1/f_a(p)$ is about $\log x.$ So actual bound could be much smaller than $O(x^{1/4}).$
This is also the conclusion drawn in Murty-Wong. Now, using $\omega(a^j-1) = O(j)$ one gets the trivial estimate $W(x) = O(x).$ If, as seems possible, the $a^j-1$ in general have a random mixture of small and large primes (rather than a mixture of predominantly small primes), then it is reasonable to expect that $W(x) = O(x^\epsilon),$ thus making $S(x)$ satisfy the Murty-Srinivasan criterion.  Indeed, in many cases one seems to have  $\omega_x(a^j-1) = 0.$

The best bound we know for $\displaystyle \sum_{p < x} \frac 1 {f_a(p)} $ is $\sqrt x /(\log x)^{1+\delta}$ where $\delta$ is some small
positive quantity, given by P. Erd\"os and M. R. Murty  \cite{Erdos-Murty}. The bounds A and B of Theorem 1, numerically at least, seem to give better estimates.

\vfill\eject

\noindent{\bf 2b. Proof of bounds A and B of Theorem 1}

\bigskip

Using summation by parts, one could reduce the estimate of $W(x)$ to that of $\displaystyle D(u)$ for $u \le x$ for positive integers $u,$ where $\displaystyle D(u) = \sum_{j \le u} \omega_x(a^j-1). $

Now it is easy to see that 

$$D(u) = \sum_{j \le u} \omega_x(a^j-1) = 
\sum_{p \le x} \frac 1 p  \sum_{j \le u} \sum_{k=1}^{p}\zeta_p^{k(a^j-1)}. \eqno{(1)}$$

\noindent {\bf Bound A:}

\medskip

Note that, for a fixed prime $p$ we have, 
on the one hand, $\displaystyle \sum_{j \le u} \omega_p(a^j-1) = \left  \lfloor\frac u {f_p} \right \rfloor .$
Here we abused the notation a little bit to make $\omega_p(m)$ equal 1 if $p $ divides $ m$ and $0$ otherwise.
This will be the case regardless of whether $f_p$ divides $u$ or not. Moreover, if $f_p > u,$ then clearly $p$ will not divide any of
the $a^j-1, j \le u,$ and $\displaystyle \left  \lfloor\frac u {f_p} \right \rfloor  = 0.$

For the part of the sum $D(u)$ with $k=p$ we have $\displaystyle \sum_{p < x} \frac 1 p  \sum_{j \le u} \sum_{k=1}^{p}\zeta_p^{k(a^j-1)}
= \sum_{p < x} \sum_{j \le u} \frac 1 p$ and this is $ O(u \log \log x).$

In what follows we focus on the $k < p$  part.  Also let $u \equiv u_f \pmod{ f_p }$ with $0 \le u_f < f_p.$ There will be no confusion owing to the dependence on $p.$

 $$D(u) = O(u \log \log x) + \sum_{p \le x} \frac 1 p \sum_{j=1}^u \sum_{k=1}^{p-1} \zeta_p^{(a^j-1)k} $$
 $$ = O(u\log \log x) + \sum_{p \le x} \frac 1 p\sum_{j=1}^{u_f}\sum_{k=1}^{p-1} \zeta_p^{(a^j-1)k} \quad \quad \qquad \qquad $$
$$\qquad \quad \quad +\sum_{p \le x} \frac 1 p \left  \lfloor\frac u {f_p} \right \rfloor \sum_{H_pg_i \in \mathbb F_p^*/H_p}
\left( \sum_{k \in H_pg_i}  \zeta_p^{-k} \left(\sum_{j=1}^{f_p} \zeta_p^{a^jk} \right)\right). \eqno{(2)}$$

In the last term of the sum in (2) above we used the facts that the sum over $j$ is periodic with period $f_p$ and also that the inner sum 
$\displaystyle \sum_{j=1}^{f_p} \zeta_p^{a^jk} $ is identical for $k$ in the same coset.

In the second term  of the sum above, for a fixed $j$ we have  $\displaystyle  \sum_{k=1}^{p-1} \zeta_p^{(a^j-1)k} = -1$ if $f_p \nmid j$ which would be the case for $j \le u_f < f_p.$ So we get
 $$D(u) 
= O(u \log \log x) + \sum_{p \le x} \left(-\frac {u_f}p \right) + 
\sum_{p \le x} \frac 1 p \left  \lfloor\frac u {f_p} \right \rfloor \sum_{H_pg_i \in \mathbb F_p^*/H_p}
\left( \sum_{k \in H_pg_i}  \zeta_p^{-k} \eta_{i,p}\right)$$
The second term of the sum is $O(u \log \log x)$ and for the last term we bound trivially and then since each coset (i.e, each $\eta_{i,p}$) appears $f_p$ times, we get the bound A:
 $$D(u) = O(u\log \log x)+ \sum_{p \le x} \frac 1 p \left \lfloor\frac u {f_p} \right \rfloor  \sum_{ i = 1} ^{e_p}f_p|\eta_{i,p}| $$
$$\implies D(u)  \le  O(u\log \log x)+ u \sum_{p \le x} \frac {\displaystyle \sum_{ i = 1} ^{e_p} |\eta_{i,p}|}p.$$
 
 \bigskip
 
\noindent {\bf Note:}  The inner sum of equation (2), namely $\displaystyle \sum_{H_pg_i \in \mathbb F_p^*/H_p}
\left( \sum_{k \in H_pg_i}  \zeta_p^{-k} \left(\sum_{j=1}^{f_p} \zeta_p^{a^jk} \right)\right),$ can be seen to equal $\sum_i \overline \eta_{i,p} \eta_{i,p} 
= \sum_i  |\eta_{i,p}|^2$ and this is known to be equal  to $p-f_p.$ See, for instance, Thaine \cite{Thaine}, p. 36. So, left as it is, this sum does not give anything new.
Although this equation actually gives a formula to compute $f_p$ in terms of the cyclotomic periods, one expects it to be hard to work with, given that $f_p$ is hard to compute even on average.

 Even though it might be bigger than original sum, the bound A (as well as bound B below) 
 at least eliminates  $f_p$ from the estimate. 
 
 \bigskip

\noindent {\bf Bound B:}

\medskip

To get bound B we use summation by parts. Since, as seen above,  the sum for $j = 1$ to $ u_f$  contributes only $O(u \log \log x)$ as does the sum for $k=p$ for each $p$, we simply focus on the remainder of the expression for $D(u).$ We call this part $D(u)'.$

$$D(u)' = \sum_{p \le x} \frac 1 p \left  \lfloor\frac u {f_p} \right \rfloor \sum_{H_pg_i \in \mathbb F_p^*/H_p}
\left( \sum_{k \in H_pg_i}  \zeta_p^{-k} \eta_{i,p}\right) $$
Before summing by parts, first we rearrange the sum so that it goes from $g_i = 1$ to $g_i= p-1.$
$$D(u)' = \sum_{p \le x} \frac 1 p \left  \lfloor\frac u {f_p} \right \rfloor \left(\sum_{g_i = 1}^{p-1} \zeta_p^{-g_i} \eta_{i,p}\right)$$

Now we rearrange the indices again and try to replace the $g_i$ with $h_m \in [1,p-1]$  as $m$ goes from $1$ to $p-1$ such that  $h_{m+1}-h_m \equiv f_p \pmod{p}.$ This can be done by letting $ h_m \equiv mf_p \pmod{p}.$  Since $p \nmid f_p$ the map $ m \to h_m$
will result in a permutation of $\{1,2,..., p-1 \}.$ 

  In particular, $\displaystyle \sum_{g_i} \eta_{i,p} = \sum_{i=1}^{p-1} \eta_{i,p} = \sum_{m=1}^{p-1} \eta_{h_m,p}.$ 
  
  Summing by parts, 
$$D(u)' = \sum_{p \le x} \frac 1 p  \left  \lfloor\frac u {f_p} \right \rfloor \left[   \left( \sum_{g_i} \eta_{i,p} \right) \zeta_p^{-h_{p-1}} 
+ \left(\sum_{k=1}^{p-2} \left( \sum_{m=1}^k \eta_{h_m,p} \right) (\zeta_p^{-h_k}-\zeta_p^{-h_{k+1}}) \right) \right]. \eqno{(3)}$$
								
Now note that $\displaystyle \left | \zeta_p^{-h_k}-\zeta_p^{-h_{k+1}} \right | = \left | 2\sin \left(\frac { \pi (h_{k+1}-h_k)}p \right) \right | \le \frac{2\pi f_p}p$ because $h_{k+1}-h_k = f_p,$ for $k = 1,2,...,p-2.$  Also note that, since $\displaystyle \sum_{i=1}^{e_p} \eta_{i,p} = -1, $ and there are $f_p$ full sets of $\eta_{i,p},$ we have $ \sum_{g_i} \eta_{i,p} = -f_p.$ Using these two computations, we get

$$\left | D(u)' \right |  \le  \sum_{p \le x} \frac 1 p \left  \lfloor\frac u {f_p} \right \rfloor \left[f_p + \left(\sum_{k=1}^{p-2} \left | \sum_{m=1}^k \eta_{h_m,p} \right | 
\frac{2\pi f_p}p \right ) \right] \le  u\sum_{p \le x} \frac 1 p  
+ 2\pi u \sum_{p \le x} \frac { \left(\sum_{k=1}^{p-2} \left | \sum_{m=1}^k \eta_{h_m,p} \right | \right) }  {p^2} $$
Adding the $O(u(\log \log x))$ terms from the sum $D(u)$ that were left out, we get bound B:
$$D(u)  \le O(u \log \log x) + 2\pi u \sum_{p \le x} \frac { \left(\sum_{k=1}^{p-2} \left | \sum_{m=1}^k \eta_{h_m,p} \right | \right) }  {p^2} $$

Due to the increased possibility of cancellations, we think bound B would be better than bound A. Numerically that seems to be the case.

\bigskip

\noindent{\bf 2c. Estimates for Bounds A and B}

\bigskip

Clearly, the best estimates for the bounds A and B would be obtained using estimates for the average values of $|\eta_{i,p}|$ and 
$\displaystyle  \left | \sum_{m=1}^k \eta_{h_m,p} \right | $ over the cyclotomic periods as well as over the primes. Here we restrict ourselves to 
seeing what estimates can be obtained for a fixed prime. As such, in this section we will drop the subscript $p$ from $\eta_{i,p}$ and 
$\eta_{h_m, p}.$

\bigskip

\noindent {\bf Bound A: } 

As mentioned in the Note in the previous section, we have $\sum_i \overline \eta_i \eta_i = \sum_i |\eta_i|^2 = p-f_p.$  From this we see that all the cyclotomic periods  have absolute values bounded by $\sqrt p, $ and possibly the absolute values are much less than $\sqrt p,$ at least on average. Indeed, we have (looking at the general term from the sum in bound A):
$$\left(\sum_i |\eta_i|\right)^2 \le e_p \sum_i |\eta_i|^2 \implies \frac 1 p \sum_i |\eta_i| \le \frac {\sqrt{e_p}}p \sqrt{\sum_i |\eta_i|^2}$$
$$\quad \quad \implies \frac 1 p \sum_i |\eta_i| \le \sqrt{ \frac {p-1}{f_p} } \frac{\sqrt{p-f_p}}p =  \sqrt{ \frac 1 {f_p} -\frac 1 p} \times \sqrt{1-\frac {f_p} p}. \eqno(4)$$
(The first inequality above follows from the elementary fact that the average of a set of positive numbers is smaller than the square root of the average of their squares).

As expected the estimate from equation (4) does not  yield a better estimate for $D(u)$ (and hence for $\sum_{p<x}1/f_p$) than the existing ones.

\medskip

On the other hand there are several estimates that show that $\displaystyle \text{max}_i \ | \eta_i| $ is smaller than $f_p$ when $f_p$ itself is small. For instance  we have $$ \text{max}_i \ | \eta_i| <  p^{ -3  \alpha /8} f_p \text{ when } f_p > p^{\frac 1 3 + \alpha}.$$ See Bourgain \cite{Bourgain} for a survey of such results. H. Montgomery, R. C. Vaughan and T. D. Wooley \cite{Montgomery} conjecture that $\displaystyle \text{max}_i \ |\eta_i| < \text{min}\left( \sqrt p, C\sqrt{f_p \times \log p} \right), $  with the constant $C$ independent of $p.$ 

Suppose $\displaystyle \text{max}_i \ |\eta_i| < M_p.$ When one uses such an estimate  in bound A one gets  
$$D(u) = \sum_{j \le u} \omega_x(a^j-1) \le  O(u\log \log x)+ u \sum_{p \le x} \frac {\displaystyle \sum_{ i = 1} ^{e_p} |\eta_{i,p}|}p$$
$$ \le O(u\log \log x)+ u \sum_{p \le x} \frac 1 p \left(\frac {p-1}{f_p}( M_p )\right) \le O(u \log \log x)+ u \sum_{p \le x} \frac {M_p}{f_p}.$$
Since one key problem is that we don't know the amount of primes with $f_p$  in a given interval $p^a \le f_p \le p^b$ it is difficult to average the bound A over the primes $p.$ Even otherwise, it appears that using existing estimates for the maximum of the periods it would be difficult to get a useful estimate for $D(u)$ and hence $\sum_{p<x} 1/{f_p}.$ One does need a good estimate for the averages of the $|\eta_i|$ as opposed to their maximum.

\bigskip

\noindent {\bf Bound B: } 

For bound B we have the following very basic estimate in the form of a trigonometric sum:

We first estimate $\displaystyle \left | \sum_{m=1}^k  \eta_{h_m,p} \right |.$ As before we omit $p$ in the subscript and denote $h_{m,p} $ by $h_m.$ We get:

$$ \quad \left | \sum_{m=1}^k \eta_{h_m} \right |  = \left | \sum_{m=1}^k \sum_{j=1}^{f_p} \zeta_p^{h_m a^j}  \right |
= \left | \sum_{m=1}^{k} \sum_{j=1}^{f_p} \zeta_p^{h_m a^j} \right |$$
$$\implies \left | \sum_{m=1}^k \eta_{h_m} \right | = \left | \sum_{j=1}^{f_p}  \sum_{m=1}^{k} \zeta_p^{mf_p a^j} \right | 
 =  \left | \sum_{j=1}^{f_p} \zeta_p^{f_pa^j} \sum_{m=0}^{k-1} \zeta_p^{mf_p a^j} \right |  
 \le \sum_{j=1}^{f_p}  \left | \frac{\sin(\pi a^jf_pk/p)}{\sin(\pi a^jf_p/p)}\right |. $$
 
Putting everything together, we get for Bound B the following estimate in terms of a trigonometric sum:

$$D(u)  \le O(u \log \log x)  + 2\pi u \sum_{p \le x} \frac 1 {p^2} \left(\sum_{k=1}^{p-2} \sum_{j=1}^{f_p}  \left | \frac{\sin(\pi a^jf_pk/p)}{\sin(\pi a^jf_p/p)}\right |\right). $$

\bigskip

\noindent{\bf 2d. Implications for Artin's Conjecture}

\bigskip

If $D(u)$ is of order $u x^{1/4}/(\log x)$ or less then after summation by parts we get that $W(x)$ is of order $ x^{1/4}$  and that would verify the Murty-Srinivasan condition and hence Artin's conjecture. 

Now, when $f_p > p^{3/4}$ the $\sum_{p<x} 1/f_p$ is at most $O(x^{1/4}/(\log x)).$  So main problem is when  $ f_p < p^{3/4}.$ 

 Erd\"os and Murty \cite{Erdos-Murty} show that for most primes $f_p$ is of order bigger than $\sqrt p.$ In the bounds given here most of the contribution would come from the primes $p$ with small $f_p$ and for such $p$ the approximation of $\sin (\pi f_p/p)$ by $\pi f_p/p$ would be better.

Finally, the full sums of bounds A and B seem to be at least as big as $ux^{1/4},$ based on numerical evidence.

Nevertheless, one could try to use the expression for  $D(u)$ in terms of cyclotomic periods (equation (2) and (3)), without taking absolute values as we did here to get the bounds A and B. If it is possible to sum over the primes first before summing over the $k$'s \`a la Voronoi in his estimate for the Dirichlet divisor problem (cf. for instance, Tenenbaum \cite{Tenenbaum} I.6.4) one could make progress on the conjecture.

In this direction would the deep estimates of exponential sums over finite fields developed by Bourgain, Konyagin and others  be useful ?
See for instance Bourgain's survey article  \cite{Bourgain}. Perhaps those methods would yield good estimates for average sums for the cyclotomic periods.

\end{document}